\RequirePackage{ifpdf}
\ifpdf 
\documentclass[pdftex]{sigma}
\else
\documentclass{sigma}
\fi

\numberwithin{equation}{section}
\numberwithin{proposition}{section}
\numberwithin{definition}{section}
\numberwithin{theorem}{section}
\numberwithin{lemma}{section}

\newcommand{\R}{\mathbb{R}}
\newcommand{\C}{\mathbb{C}}

\newcommand{\N}{\mathbb{N}}

\newcommand\dint{\displaystyle\int}

\begin{document}

\renewcommand{\thefootnote}{$\star$}

\renewcommand{\PaperNumber}{084}

\FirstPageHeading

\ShortArticleName{Dunkl  Hyperbolic Equations}
\ArticleName{Dunkl  Hyperbolic Equations\footnote{This paper is a contribution to the Special
Issue on Dunkl Operators and Related Topics. The full collection
is available at
\href{http://www.emis.de/journals/SIGMA/Dunkl_operators.html}{http://www.emis.de/journals/SIGMA/Dunkl\_{}operators.html}}}
\Author{Hatem MEJJAOLI}

\AuthorNameForHeading{H. Mejjaoli}

\Address{Faculty of Sciences of Tunis, Department of Mathematics,
1060 Tunis, Tunisia}

\Email{\href{mailto:hatem.mejjaoli@ipest.rnu.tn}{hatem.mejjaoli@ipest.rnu.tn}}

\ArticleDates{Received May 10, 2008, in f\/inal form November 24,
2008; Published online December 11, 2008}

\Abstract{We introduce and study  the Dunkl  symmetric systems.
 We prove the well-posedness results for the Cauchy problem for these  systems.
 Eventually we describe the f\/inite speed of it.
  Next the semi-linear Dunkl-wave equations are also studied.}


\Classification{35L05;  22E30}

\section{Introduction}

We consider the dif\/ferential-dif\/ference
operators $T_j$, $j = 1,\! \dots ,\! d$, on $\R^{d}$ introduced by
C.F.~Dunkl in~\cite{D1} and called Dunkl operators in the
literature.  These
operators are very important in pure mathematics and in physics.
They  provide a useful tool in the study of special functions with
root systems~\cite{D2}.

 In
this paper, we are interested  in studying two types of Dunkl
hyperbolic equations. The f\/irst one is the Dunkl-linear symmetric system
\begin{gather}
\left\{
\begin{array}{l}
\displaystyle  \partial_t u - \sum_{j=1}^{d}A_j T_j u - A_0 u  =  f, \\
  u|_{t=0}  =  v,
\end{array}
\right.
\label{S}
\end{gather} where the $A_p$ are square matrices  $m\times m$ which satisfy some hypotheses (see Section~\ref{sec3}),
the initial data  belong to Dunkl--Sobolev spaces $[H^{s}_{k}(\R^d)]^m$ (see~\cite{MT1}) and $f$ is a continuous function on an interval $I$ with value in $[H^{s}_{k}(\R^d)]^m$. In the classical case  the Cauchy problem for symmetric hyperbolic systems of f\/irst order has been introduced and studied by Friedrichs~\cite{F1}. The Cauchy problem will be solved with the aid of energy integral inequalities, developed for this purpose by Friedrichs. Such energy inequalities have been employed by H.~Weber~\cite{W}, Hadamard~\cite{Ha}, Zaremba~\cite{Z}  to derive various uniqueness theorems, and by Courant--Friedrichs--Lewy~\cite{CFL},   Friedrichs~\cite{F1}, Schauder~\cite{Sc} to derive existence theorems. In all these treatments the energy inequality is used to show that the solution, at some later time, depends boundedly  on the initial values in an appropriate norm. However, to derive an existence theorem one needs, in addition to the a priori energy estimates, some auxiliary constructions. Thus, motivated by these methods
 we will  prove by energy methods  and Friedrichs approach  local
well-posedness and principle of f\/inite speed of propagation for the system \eqref{S}.

Let us f\/irst summarize our well-posedness results and  f\/inite speed of propagation (Theorems~\ref{theorem3.1} and~\ref{theorem3.2}).

{\bf{Well-posedness for DLS.}} For all given  $f \in
[C(I,H^{s}_{k}(\R^d))]^m$ and $v \in [H^{s}_{k}(\R^d)]^m$,
 there exists a unique solution $u$ of the system \eqref{S} in the
 space
\[
 [C(I,H^{s}_{k}(\R^d))]^m\bigcap
[C^{1}(I,H^{s-1}_{k}(\R^d))]^m.
\]
In the classical case,  a similar   result can be found in~\cite{Chaz}, where the authors used another method based on the symbolic calculations  for the pseudo-dif\/ferential operators that we cannot adapt for the system \eqref{S} at the moment. Our method uses some ideas inspired by the works~\cite{Chaz,F,Fr,F1,Fri,FL1,FL2,K,Lax,LPh,Ma,Rau}. We note that K.~Friedrichs has solved the Cauchy problem in a~lens-shaped domain \cite{F1}. He proved existence of extended solutions by Hilbert space method and showed the dif\/ferentiability of these solutions using  complicated calculations.
A~similar problem is that of a symmetric hyperbolic system studied by P. Lax, who gives a method  of\/fering  both the existence and the dif\/ferentiability of solutions at once~\cite{Lax}.
He reduced the problem to the case where all functions are periodic in every independent variable.

{\bf{Finite speed of propagation.}} Let \eqref{S} be as above. We
assume that  $f \in [C(I, L^{2}_k(\R^d))]^m$ and $v
\in [L^{2}_k(\R^d)]^m$.

 $\bullet$ There exists a positive constant $C_0 $  such that, for any
positive real  $R$   satisfying
\[
\left\{\begin{array}{ll}\label{10} f(t,x)\equiv 0\
&\mbox{for} \ \  \|x\| < R - C_0 t,\\  v(x) \equiv 0\  &\mbox{for} \ \
 \|x\| < R,
\end{array}
\right.
\] the unique solution $u$ of the system  \eqref{S}  satisf\/ies
\[
u(t,x) \equiv 0 \qquad \mbox{for} \quad \|x\| < R - C_0 t.
\]

$\bullet$ If given $f$ and $v$ are such that
\[
\left\{\begin{array}{ll}\label{10} f(t,x)\equiv 0\
&\mbox{for} \ \  \|x\| > R + C_0 t,\\  v(x) \equiv 0\  &\mbox{for} \ \
\|x\| > R,
\end{array}
\right.
\] then the unique solution $u$ of the system  \eqref{S}  satisf\/ies
\[
u(t,x)\equiv 0 \qquad \mbox{for} \quad \|x\| > R + C_0 t.
\]
In the classical case,  similar  results can be found in~\cite{Chaz} (see also~\cite{Shi}).

 A standard
example of the Dunkl linear symmetric system is the
Dunkl-wave equations with variable coef\/f\/icients def\/ined
by
\[
   \partial^{2}_{t}u - {\rm{div}}_{k}[A\cdot \nabla_{k,x }u] +
  Q(t,x,\partial_t u, T_x u), \qquad t \in \R, \qquad x \in \R^d,
\]
where
\[
\nabla_{k,x}u = \left(T_1\, u,\dots,
T_d\, u\right), \qquad {\rm{div}}_{k}  \left(  v_1,\dots,v_d  \right) =
 \sum_{i=1}^{d}T_i v_i,
\]  $A$ is a real
symmetric
 matrix which satisf\/ies some hypotheses (see Subsection~\ref{sec3.2})
and $Q(t,x,\partial_tu$, $T_x u)$ is dif\/ferential-dif\/ference operator of degree~1 such that these  coef\/f\/icients are $C^{\infty}$, and all derivatives are bounded.

      From the previous results we deduce the
well-posedness of the generalized Dunkl-wave equations (Theorem~3.3).

{\bf{Well-posedness for GDW.}} For all $s  \in \N$, $u_0 \in
H^{s+1}_{k}(\R^d)$, $u_1 \in H^{s}_k(\R^d)$ and $f$  in $C(\R,
H^{s}_k(\R^d))$, there exists a unique $u \in
C^{1}(\R,H^{s}_k(\R^d)) \cap C(\R,H^{s+1}_k(\R^d))$ such that
\[
\left\{
\begin{array}{l}
  \partial^{2}_{t}u - {\rm{div}}_{k}[A\cdot \nabla_{k,x }u] +
  Q(t,x,\partial_t u, T_x u)   =   f, \\
  u|_{t=0}   =   u_0,\\
\partial_t u|_{t=0}   =   u_1.
\end{array}
\right.
\]

The second type of Dunkl hyperbolic equations that we are interested
is the  semi-linear Dunkl-wave equation
\begin{gather}
\left\{
\begin{array}{l}
  \partial_{t}^{2} u - \triangle_k u   =   Q(\Lambda_k u , \Lambda_k u), \\
  (u,\partial_{t} u)|_{t=0}   =   (u_{0},u_1),
\end{array}
\right.
\label{P}
\end{gather}
 where
 \[
 \triangle_k = \sum_{j = 1}^d T_j^2,  \qquad \Lambda_k u =
(\partial_{t} u, T_1 u, \dots, T_d u),
\]
 and $Q$ is a quadratic form on $\R^{d+1}$.

Our  main result for this type of Dunkl hyperbolic equations
is the following.

 {\bf{Well-posedness for SLDW.}} Let $(u_0,u_1)$ be in
$ H^{s}_{k}(\R^d)\times H^{s-1}_{k}(\R^d)$ for $s > \gamma +
\frac{d}{2} + 1$. Then there exists a positive time $T$  such that the problem \eqref{P} has a unique solution $u$ belonging to
\[
  C([0,T],H^{s}_{k}(\R^d))\cap
C^{1}([0,T],H^{s-1}_{k}(\R^d))
\]  and satisfying  the blow up criteria  (Theorem~\ref{theorem4.1}).

In the classical case see \cite{B1,B2,B3,S1}.
We note that the Huygens' problem for  the homogeneous Dunkl-wave equation is studied by S. Ben Sa\"{i}d and B. {\O}rsted~\cite{Ben}.

 The
 paper is organized as follows.
In Section~\ref{sec2}  we recall  the main results about the
harmonic analysis associated with the Dunkl operators. We study in Section~\ref{sec3}  the
generalized Cauchy
 problem of the Dunkl linear  symmetric systems, and we
prove the  principle of  f\/inite speed of propagation of these
systems.
In the last section we study  a semi-linear Dunkl-wave
 equation and we prove the well-posedness of this equation.

 Throughout this paper by $C$
we always represent a positive constant not
      necessarily the same in each occurrence.

\section{Preliminaries}\label{sec2}

This section gives an introduction to the theory of Dunkl operators, Dunkl transform, Dunkl convolution and to the Dunkl--Sobolev spaces. Main references are \cite{D1,D2,D3,J,MT1,MT2,R1,Ros,karni4,T1,T5}.

\subsection[Reflection groups, root systems and multiplicity functions]{Ref\/lection groups, root systems and multiplicity functions}\label{sec2.1}

The basic ingredient in the theory of Dunkl operators are root systems and f\/inite ref\/lection groups, acting on  ${\R}^{d}$ with the standard Euclidean
scalar product $\langle\cdot,\cdot\rangle$ and $||x|| = \sqrt{\langle x , x
\rangle}$. On~${\C}^{d}$,  $\|\cdot\|$ denotes also the standard Hermitian norm, while $\langle z,w\rangle= \sum\limits_{j=1}^{d}z_{j}\overline{w}_{j}$.

For $\alpha \in
\R^{d}\backslash\{0\} $, let $\sigma_{\alpha}$ be the ref\/lection
in the
 hyperplane  $H_{\alpha} \subset \R^{d} $ orthogonal to $\alpha$, i.e.
\begin{equation*}
\sigma_{\alpha}(x) = x - 2 \frac{\langle\alpha, x\rangle} {||\alpha||^{2}}
\alpha. 
\end{equation*}

 A f\/inite set  $R \subset \R^{d}\backslash\{0\} $
 is called a root system if
$R \cap {\R}\cdot\alpha = \{\alpha, - \alpha\}$ and
$\sigma_{\alpha} R = R$ for all $\alpha \in R$. For a given root
system $R$
 the ref\/lections $\sigma_{\alpha}$, $\alpha \in R $,
generate a f\/inite group $W \subset O(d)$, called the ref\/lection
group associated with $R$. All ref\/lections in $W$ correspond to suitable pairs of roots.
We f\/ix a positive root system
$R_{+} = \big\{\alpha \in R\, / \langle\alpha, \beta\rangle > 0 \big\}$ for some  $\beta \in \R^{d} \backslash
\displaystyle\bigcup_{\alpha \in R}H_{\alpha}$. We will assume that $\langle\alpha,\alpha\rangle =
2$ for all $\alpha \in R_{+}$.

 A function $k :R \longrightarrow \C$  is called a multiplicity function if it is invariant
under the action of the associated ref\/lection group $W$.
  For abbreviation, we introduce the index
\begin{equation*}
\gamma = \gamma(k) =  \sum_{\alpha \in R_{+}
}k(\alpha). 
\end{equation*}
Throughout this paper, we will assume
 that the multiplicity is non-negative, that is $k(\alpha) \geq 0$ for all $\alpha \in R$. We write $k \geq 0$ for short.
Moreover, let $\omega_{k}$ denote the weight function
 \begin{equation*}
\omega_{k}(x) = \prod_{\alpha \in R_{+}
}|\langle\alpha,x\rangle|^{2k(\alpha)}, 
\end{equation*}
which is invariant and homogeneous of degree $2 \gamma$. We introduce the Mehta-type constant
\begin{equation*}
c_{k} = \left(\int_{\R^{d}}\exp(-||x||^{2})\omega_{k}(x)\,dx\right)^{-1}.
\end{equation*}

\subsection{The Dunkl operators and the Dunkl kernel}\label{sec2.2}

  We denote by
 \begin{enumerate}\itemsep=0pt
\item[--] $C(\R^{d})$  the space of continuous functions on $\R^{d}$;
\item[--] $C^{p}(\R^{d})$ the space of functions of class $C^p$ on $\R^{d}$;
\item[--] $C_{b}^{p}(\R^d)$ the space of bounded functions of class $C^p$;
\item[--] $ {\cal E}(\R^{d})$ the space of
$C^{\infty}$-functions on $\R^{d}$;
\item[--] ${\cal S}(\R^{d})$   the Schwartz
space of  rapidly
decreasing functions on $\R^{d}$;
\item[--] $D(\R^{d})$ the space of
$C^{\infty}$-functions  on $\R^{d}$ which are of compact
 support;
\item[--] ${\cal S'}(\R^{d})$
the space of temperate distributions on $\R^{d}$. It is the
topological dual of ${\cal S}(\R^{d})$.
\end{enumerate}

In this subsection we collect some notations and results on the Dunkl
operators
 (see \cite{D1,D2} and~\cite{D3}).
The Dunkl operators $T_{j}$, $j = 1, \dots,  d $, on $\R^{d}$
associated with the
 f\/inite ref\/lection group~$W$ and multiplicity function $k$ are given by
\begin{equation*}
T_{j} f(x) = \frac{\partial f}{\partial x_{j}}(x) +
\displaystyle\sum_{\alpha \in R_{+}}k(\alpha) \alpha_{j}
\frac{f(x) - f(\sigma_{\alpha}(x))}{\langle\alpha,x\rangle},\qquad f   \in
C^{1}(\R^{d}). 
\end{equation*}
 Some properties of the $T_j,$ $j = 1,\dots, d$,
are given in the following:

For all $f$ and $g$ in $C^{1}(\R^{d})$ with at least one of them  is $W$-invariant, we have
\begin{equation}\label{ki}
T_{j}(fg) = (T_{j}f)g + f(T_{j}g),\qquad j = 1, \dots, d.
\end{equation}

 For $f $ in   $ C_{b}^{1}(\R^{d})$ and $g$ in $
{\cal S}(\R^{d})$ we have
\begin{equation}\label{kii}
\int_{\R^{d}} T_{j}f(x) g(x)\omega_{k}(x)\,dx = - \int_{\R^{d}}
 f(x) T_{j}g(x)\omega_{k}(x)\,dx, \qquad j = 1, \dots, d.
\end{equation}

We def\/ine the Dunkl--Laplace operator on $\R^{d}$ by
\begin{equation*}
\triangle_{k}f(x)  =  \sum
_{j = 1}^{d}T_j^2 f(x) = \triangle f(x) + 2
 \sum_{\alpha \in R^+}
k(\alpha)\left[\dfrac{\langle\nabla
f(x),\alpha\rangle}{\langle\alpha,x\rangle} - \dfrac{f(x) -
f(\sigma_{\alpha}(x))}{\langle\alpha,x\rangle^2}\right]. 
\end{equation*}

For $y \in \R^{d} $,  the system \begin{equation*} \left\{
\begin{array}{l}
T_{j}u(x,y) = y_{j} u(x,y), \qquad j = 1, \dots, d,\\ u(0,y) = 1,
\end{array}
\right.
\end{equation*} admits a unique analytic solution on $\R^{d}$, which
will be denoted by $K(x,y)$ and called Dunkl kernel. This kernel has
a unique holomorphic extension to ${\C}^{d} \times {\C}^{d}$.

The Dunkl kernel possesses the following properties:

i) For  $z, t \in \C^d$, we have $K(z,t) = K(t,z);$
$K(z,0) = 1$ and $K(\lambda z, t) = K(z, \lambda t)$ for all
$\lambda \in \C$.

ii) For all $\nu \in \N^d$, $x \in
{\R}^d$ and $z \in \C^d$ we have
\[
|D^\nu_z K(x,z)| \leq
||x||^{|\nu|} \exp(||x|| \, ||{\rm{Re}}\,z||),
\] with
\[
D^\nu_z =
\displaystyle{\frac{\partial^{|\nu|}}{\partial z_1^{\nu_1}\cdots
\partial z_d^{\nu_d} }}\qquad {\rm{and}} \qquad |\nu| = \nu_1 + \cdots + \nu_d.
\]
 In
particular for all $x, y \in {\R}^d$:
\[
|K(-ix, y)| \leq 1.
\]

iii) The function $K(x,z)$ admits for all $x \in
{\R}^d$ and $z \in \C^d$ the following Laplace type integral
representation \begin{equation}K(x,z) = \int_{{\mathbb R}^d} e^{\langle y,
z\rangle}d\mu_x(y),\label{6}
\end{equation} where $\mu_x$ is a probability
measure on ${\R}^d$ with support in the closed ball $B(0, ||x||)$
of center $0$ and radius $\|x\|$ (see~\cite{R1}).

The Dunkl intertwining operator $V_k$ is the operator from
   $C({\R}^d)$ into itself given by
   \begin{equation*}
 V_k f(x) = \int_{\R^d}
f(y)d\mu_{x}(y), \qquad {\rm{for \  all}} \ \ x \in \R^{d},  
\end{equation*}
   where $\mu_x$ is the measure given by the  relation (\ref{6}) (see \cite{R1}). In
   particular, we have
\[
 K(x,z) =
   V(e^{\langle \cdot, z\rangle})(x), \qquad {\rm{for \  all}}\ \ x \in {\R}^d \quad {\rm{and}}\ \ z \in \C^d.
   \]
In \cite{D2} C.F.~Dunkl proved that
 $V_k$  is a linear isomorphism from the space of homogeneous polynomial
${\cal P}_n$ on $\R^{d}$ of degree n into itself satisfying the relations
\begin{gather}\label{H1}
\left\{
\begin{array}{l}
T_j V_k  =  V_k \dfrac{\partial}{\partial x_j}, \qquad j = 1, \dots , d, \vspace{1mm}\\
V_k(1)  =  1.
\end{array}
\right.
\end{gather}

K.~Trim\`eche has proved in  \cite{T1} that the operator $V_k$ can be extended to
 a topological isomorphism from ${\cal E}(\R^{d})$  into itself satisfying the relations~(\ref{H1}).

 \subsection{The Dunkl transform}\label{sec2.3}

 We denote by
$L_{k}^{p}(\R^{d})$ the space of measurable functions on $\R^{d}$
such that
\begin{gather*}
||f||_{L_{k}^{p}(\R^{d})}  :=  \left(  \int_{\R^{d}} |f(x)|^{p}
\omega_{k}(x) \,dx\right)^{\frac{1}{p}} < +\infty
 \qquad \mbox{if} \quad  1 \leq p
< + \infty,\\
||f||_{L_{k}^{\infty}(\R^{d})}  :=   {\rm{ess}}\;  \sup _{x \in \R^{d} }
|f(x)| < +\infty.
\end{gather*}

The Dunkl transform of a function $f$ in $L_{k}^{1}(\R^{d})$ is
given by
\begin{equation*}
 {\cal F}_{D}(f) (y) =
\int_{\R^{d}}f(x) K(-iy,x) \omega_{k}(x)dx,\qquad {\rm{for \  all}} \ \  y \in \R^{d}. 
\end{equation*}
 In the following we give some  properties of this transform (see
\cite{D3,J}).

 i) For $f$
in $L_{k}^{1}(\R^{d})$ we have
\begin{equation*}
||{\cal F}_{D} (f)|| _{L_{k}^{\infty}(\R^{d})} \leq ||f||_{L_{k}^{1}(\R^{d})}. 
\end{equation*}

 ii) For $f$ in ${\cal S}(\R^{d})$ we have
\begin{equation*}
  {\cal F}_{D}( T_j f)( y ) = i y_j
{\cal F}_{D}( f ) (y), \qquad {\rm{for \  all}} \ \   j = 1, \dots,d \quad {\rm{and}} \quad y \in \R^{d}. 
\end{equation*}

\subsection{The Dunkl convolution}\label{sec2.4}

\begin{definition} \label{definition2.1} Let $y$ be in $\R^{d}$. The Dunkl translation operator $f \mapsto \tau_y f$ is
def\/ined on ${\cal S}(\R^d)$ by
\begin{gather*}
{\cal F}_D (\tau_{y}f)(x)= K(ix,y){\cal F}_D (f)(x), \qquad {\rm{for \  all}} \ \ x \in \R^{d}. %
\end{gather*}
\end{definition}

\begin{proposition}\label{proposition2.1} {\rm{i)}}
  The  operator $\tau_y$, $y
\in \R^d$, can also be  def\/ined    on ${\cal E}(\R^d)$ by
\begin{equation*}
\tau_{y}f(x)= (V_k)_x (V_k)_y[(V_k)^{-1}(f)(x+y)], \qquad {\rm{for \  all}} \ \  x \in \R^{d}
\end{equation*}
(see {\rm \cite{T5}}).

 {\rm{ii)}} If $f(x) = F(||x||)$ in ${\cal E}(\R^{d})$, then we have
 \[
 \tau_{y}f(x) = V_k \left[F
(\sqrt{||x||^2 + ||y||^2 +2 \langle x,\cdot\rangle
})\right](x), \qquad {\rm{for \  all}} \ \ x \in \R^{d}
\]
 (see {\rm \cite{Ros}}).
 \end{proposition}

Using the Dunkl
translation operator, we  def\/ine the Dunkl convolution product of
functions as follows (see \cite{karni4,T5}).
\begin{definition} \label{definition2.2}
The
Dunkl convolution product of  $f$ and $g$ in
${\cal S}(\R^{d}) $ is the function $f*_{D}g$ def\/ined by
\begin{equation*}
f*_{D}g(x) = \int_{\R^{d}}\tau_{x}f(-y)g(y)\omega_{k}(y)dy, \qquad {\rm{for \  all}} \ \ x \in \R^{d}.
\end{equation*}
\end{definition}

\begin{definition}  \label{definition2.3} The Dunkl transform of a
distribution $\tau$ in ${\cal S}'(\R^{d})$ is def\/ined by
\begin{equation*}
 \langle{\cal F}_{D}(\tau), \phi \rangle = \langle \tau,{\cal F}_{D }(\phi)\rangle, \qquad {\rm{for \  all}} \ \
\phi \in {\cal S}(\R^{d}). \end{equation*}
\end{definition}

\begin{theorem}\label{theorem2.1}
 The Dunkl transform ${\cal F}_{D}$ is a topological isomorphism
from ${\cal S'}(\R^{d})$ onto itself.
 \end{theorem}

\subsection[The Dunkl-Sobolev spaces]{The Dunkl--Sobolev spaces}\label{sec2.5}

In this subsection we
recall some def\/initions  and results on Dunkl--Sobolev spaces (see
\cite{MT1,MT2}).

Let $\tau$ be in ${\cal S'}(\R^d)$. We def\/ine the
distributions $T_j \tau$, $j=1,\dots,d,$ by
\begin{gather*}
\langle T_j \tau, \psi\rangle   =   -
\langle\tau, T_j \psi\rangle, \qquad \mbox{for all} \ \  \psi   \in   {\cal
S}(\R^d).\end{gather*}
These distributions  satisfy the following property
\begin{gather*}
{\cal F}_{D}(T_j \tau)  =  i y_j {\cal F}_{D}( \tau),\qquad j = 1,
\dots, d. 
\end{gather*}
\begin{definition}
Let $s \in \R$, we def\/ine the space $H^{s}_{k}(\R^d)$  as the set of distributions $
u \in {\cal S'}(\R^{d})$ satisfying $ (1 + ||\xi||^2)^{\frac{s}{2}} {\cal
F}_{D}( u) \in L_{k}^{2}(\R^d).$
\end{definition}

We provide this space with
the scalar product
\begin{equation*}
\langle u,v\rangle_{H^{s}_{k}(\R^d)} = \dint_{\R^d}( 1 + ||\xi||^2)^s  {\cal
F}_{D}( u)(\xi) \overline{{\cal F}_{D}(v
)(\xi)}\omega_{k}(\xi)d\xi 
\end{equation*}
and the norm
\begin{equation*}
||u||_{H^{s}_{k}(\R^d)}^2 = \langle u,u\rangle_{H^{s}_{k}(\R^d)}.  
\end{equation*}

\begin{proposition}\label{proposition2.2} {\rm{i)}} For $s \in \R$ and $\mu \in \N^d$, the Dunkl operator
$T^{\mu}$ is continuous from $H^{s}_{k}(\R^d)$ into
$H^{s-|\mu|}_{k}(\R^d)$.

{\rm{ii)}} Let $p \in \N$. An
element $u$ is in $H^{s}_{k}(\R^d)$ if and only if for all $\mu
\in \N^d$, with $|\mu| \leq p$, $T^{\mu} u$ belongs to
$H^{s-p}_{k}(\R^d)$, and we have
\[
||u||_{H^{s}_{k}(\R^d)} \sim \sum_{|\mu|
\leq p} ||T^{\mu} u||_{H^{s-p}_{k}(\R^d)}.
\]
\end{proposition}

\begin{theorem}\label{theorem2.2} {\rm{i)}} Let $u$ and $ v \in H^{s}_{k}(\R^d)\bigcap
L^{\infty}_{k}(\R^d) $, $s > 0$, then $uv \in H^{s}_{k}(\R^d)$ and
\[
||uv||_{H^{s}_{k}(\R^d)} \leq C(k,s)
\big[||u||_{L^{\infty}_{k}(\R^d)}||v||_{H^{s}_{k}(\R^d)} +
||v||_{L^{\infty}_{k}(\R^d)}||u||_{H^{s}_{k}(\R^d)}\big].
\]

{\rm{ii)}} For $s > \frac{d}{2} + \gamma$,
$H^{s}_{k}(\R^d)$ is an algebra with respect to pointwise multiplications.
\end{theorem}

\section{Dunkl linear symmetric systems}\label{sec3}

For any interval $I$ of $\R$  we def\/ine the mixed space-time spaces
 $C(I,H^s_{k}(\R^d))$, for $s \in \R$, as
 the spaces of  functions from $I$
into $H^s_{k}(\R^d)$ such that the map
\[
t \mapsto ||u(t,\cdot )||_{H^s_{k}(\R^d)}
\] is continuous.  In this section, $I$ designates the interval
$[0,T[$, $T > 0$ and
\[
u = (u_1,\dots,u_m), \qquad u_p \in
C(I,H^{s}_{k}(\R^d)),
\]
 a vector with $m$ components elements of $
C(I,H^{s}_{k}(\R^d))$. Let $(A_p)_{0 \leq p \leq d}$ be a family of
 functions from $I\times\R^d$ into the space of  $m\times m$
matrices with real coef\/f\/icients $a_{p,i,j}(t,x)$ which are $W$-invariant
with  respect to $x$ and whose all derivatives in $x \in
\R^d$ are bounded and continuous functions of~$(t,x)$.

 For a given $f \in [C(I,H^{s}_{k}(\R^d))]^m$ and $v \in
[H^{s}_{k}(\R^d)]^m$, we f\/ind $u \in
[C(I,H^{s}_{k}(\R^d))]^m$ satisfying the system  \eqref{S}. 

 We shall f\/irst def\/ine the notion of symmetric systems.
\begin{definition} \label{definition3.1} The system \eqref{S} is symmetric, if and only if, for any
$ p \in \{1,\dots,d\}$ and  any $(t,x) \in I \times \R^d$ the
matrices $A_p (t,x)$ are symmetric, i.e.\ $a_{p,i,j}(t,x)=
a_{p,j,i}(t,x)$.
\end{definition}

In this section, we shall assume  $s \in \N$ and
denote  by $||u(t)||_{s,k}$ the norm def\/ined by
\[
||u(t)||_{s,k}^{2} = \sum_{\substack{1 \leq p\leq
m\\ 1\leq |\mu| \leq
s}}||T_{x}^{\mu}u_{p}(t)||_{L^{2}_{k}(\R^d)}^{2}.
\]

\subsection{Solvability for Dunkl linear symmetric systems}\label{sec3.1}

 The aim
 of this subsection is to prove the
following theorem.
\begin{theorem}\label{theorem3.1}
  Let
\eqref{S} be a symmetric system. Assume that $f$  in $[C(I,H^{s}_{k}(\R^d))]^m$ and $v$ in  $[H^{s}_{k}(\R^d)]^m$,
then there exists a unique solution $u$ of  \eqref{S} in
$ [C(I,H^{s}_{k}(\R^d))]^m\!\!\bigcap [C^{1}\!(I,H^{s-1}_{k}(\R^d))]^m.\!$
\end{theorem}

The proof of this theorem will be made in several steps:
\begin{enumerate}\itemsep=0pt

\item[{\bf{A.}}] We prove a priori estimates  for the regular solutions
of the system \eqref{S}.

\item[{\bf{B.}}]  We apply the Friedrichs method.

 \item[{\bf{C.}}]  We pass to the limit for regular solutions
and we obtain the existence  in all cases by the
regularization of the Cauchy data.

\item[{\bf{D.}}]
 We prove the uniqueness  using the
 existence result on the adjoint system.
\end{enumerate}

{\bf{A. Energy
estimates.}}
The symmetric hypothesis is crucial for the energy
 estimates which are  only valid  for regular solutions. More precisely we have:

\begin{lemma} \label{lemma3.1} {\rm (Energy
Estimate in $[H_{s}^k(\R^d)]^{m}$).} For any positive integer $s$, there
exists a positive constant $\lambda_s$  such that, for any function $u$
in $[C^{1}(I,H^{s}_{k}(\R^d))]^m \bigcap
[C(I,H^{s+1}_{k}(\R^d))]^m$, we have
\begin{equation}||u(t)||_{s,k} \leq e^{\lambda_s
t} ||u(0)||_{s,k}  + \dint_{0}^{t}
e^{\lambda_s
(t-t')}||f(t')||_{s,k}dt', \qquad {\rm{ for \ all}} \ \  t \in I,\label{1+}
\end{equation}
with
\[
f =
\partial_t u - \sum_{p=1}^{d}A_p T_p u - A_0 u.
\]
\end{lemma}

To prove Lemma~\ref{lemma3.1}, we need the following lemma.

\begin{lemma}\label{lemma3.2} Let
$g$ be a $C^1$-function on $[0,T[$, $a$ and $b$ two positive
continuous functions.  We assume
\begin{equation}\label{2+}
\frac{d}{dt} \,g^2(t) \leq 2
  a(t)g^2(t)+ 2b(t) |g(t)|.
\end{equation}
Then, for $t \in [0,T[$, we have
\[
|g(t)| \leq |g(0)|\exp \int^t_0 a(s)ds + \int^t_0 b(s) \exp\left(\int^t_s
 a(\tau)d\tau\right)ds.
 \]
\end{lemma}

\begin{proof} To prove this lemma, let us set for $\varepsilon > 0$, $g_{\varepsilon}(t) = \big(g^2(t) + \varepsilon\big)^{\frac{1}{2}}$; the function $g_\varepsilon$ is $C^{1}$, and we have $|g(t)| \leq g_{\varepsilon}(t)$. Thanks to the inequality (\ref{2+}), we have
\[
\frac{d}{dt}(g^{2})(t)\leq  2a(t)g_{\varepsilon}^2(t)+ 2b(t) g_{\varepsilon}(t).
\]
As $\frac{d}{dt}(g^{2})(t)=\frac{d}{dt}(g_{\varepsilon}^{2})(t)$. Then
\[
\frac{d}{dt}(g_{\varepsilon}^{2})(t) = 2g_{\varepsilon}(t) \frac{dg_{\varepsilon}}{dt}(t)\leq  2a(t)g_{\varepsilon}^2(t)+ 2b(t) g_{\varepsilon}(t).
\]
Since for all $t \in [0,T[$ $g_{\varepsilon}(t) > 0$,  we deduce then
\[
  \frac{dg_{\varepsilon}}{dt}(t)\leq  a(t)g_{\varepsilon}(t)+ b(t).
\]
Thus
\[
\frac{d}{dt}\left[g_{\varepsilon}(t)  \exp\left(-\dint_{0}^{t}a(s)ds\right)\right] \leq b(t)\exp\left(-\dint_{0}^{t}a(s)ds\right).
\]
So, for $t \in [0,T[$,
\[
g_{\varepsilon}(t) \leq g_{\varepsilon}(0)\exp \dint_{0}^{t}a(s)ds + \dint_{0}^{t}b(s)\exp \left(\dint_{s}^{t}a(\tau)d\tau\right)ds.
\]
Thus, we obtain the conclusion of the lemma by tending $\varepsilon$ to zero.
\end{proof}

\begin{proof}[Proof of Lemma~\ref{lemma3.1}.]  We prove this estimate by
induction on $s$. We f\/irstly assume  that $u$ belongs to
$[C^{1}(I,L^{2}_{k}(\R^d))]^m \bigcap [C(I,H^{1}_{k}(\R^d))]^m$.
We then have  $f \in [C(I,L^{2}_{k}(\R^d))]^m$, and the
function  $t \mapsto ||u(t)||_{0,k}^2$ is  $C^1$ on the interval $I$.

By def\/inition of $f$ we have
\begin{gather*}
\frac{d}{dt}||u(t)||_{0,k}^2  =  2\langle\partial_t u,
u\rangle_{L^{2}_{k}(\R^d)} =  2\langle f,u\rangle_{L^{2}_{k}(\R^d)} +
2\langle A_0 u, u\rangle_{L^{2}_{k}(\R^d)} +
2\displaystyle\sum_{p=1}^{d}\langle A_p T_p u,
u\rangle_{L^{2}_{k}(\R^d)}.
\end{gather*} We will estimate the third term of the sum above
by using the symmetric hypothesis  of the matrix $A_p$. In fact from (\ref{ki}) and (\ref{kii}) we
have
\begin{gather*}
  \langle A_p T_p u, u\rangle_{L^{2}_{k}(\R^d)}   =
  \sum_{1\leq i,j\leq m}\dint_{\R^d}a_{p,i,j}(t,x)
  [(T_p)_x\, u_{j}(t,x)] u_{i}(t,x)\omega_{k}(x)dx \\
\phantom{\langle A_p T_p u, u\rangle_{L^{2}_{k}(\R^d)}}{}
=   -  \sum_{1\leq i,j\leq m}\dint_{\R^d}a_{p,i,j}(t,x)
 [(T_p)_x\, u_{i}(t,x)] u_{j}(t,x)\omega_{k}(x)dx\\
\phantom{\langle A_p T_p u, u\rangle_{L^{2}_{k}(\R^d)}=}{}  - \sum_{1\leq i,j\leq m}\dint_{\R^d}[(T_p)_x\,
   a_{p,i,j}(t,x)]
   u_{j}(t,x) u_{i}(t,x)\omega_{k}(x)dx.
\end{gather*}
The matrix $A_p$ being symmetric, we have
\[
- \sum_{1\leq i,j\leq
m}\dint_{\R^d}a_{p,i,j}(t,x)
  T_p u_{i}(t,x) u_{j}(t,x)\omega_{k}(x)dx= -\langle A_p T_p u, u\rangle_{L^{2}_{k}(\R^d)}.
\]
  Thus
  \[
   \langle A_p T_p u, u\rangle_{L^{2}_{k}(\R^d)} = -\frac{1}{2} \sum_{1\leq
i,j\leq m}\dint_{\R^d}T_p\, a_{p,i,j}(t,x)
   u_{i}(t,x) u_{j}(t,x)\omega_{k}(x)dx.
\]
   Since the coef\/f\/icients
   of the matrix $A_p$,  as well as their derivatives are bounded on $\R^d$ and
   continuous on $I  \times \R^d$, there exists  a positive constant $\lambda_0$
    such that
\begin{equation*}
  \dfrac{d}{dt}||u(t)||_{0,k}^{2} \leq 2
  ||f(t)||_{0,k}||u(t)||_{0,k}+ 2\lambda_0 ||u(t)||_{0,k}^{2}.
\end{equation*}
To complete the   proof of  Lemma \ref{lemma3.1} in the case $s = 0$ it
suf\/f\/ices to apply   Lemma~\ref{lemma3.2}.
 We assume now that  Lemma~\ref{lemma3.1} is proved for  $s$.

 Let $u$ be the function of  $[C^1(I, H_k^{s+1}(\R^d))]^m \cap [C(I, H_k^{s+2}(
 \R^d))]^m$, we  now introduce the function (with $m(d+1)$
 components) $U$ def\/ined by
\[
U = (u,T_1 u,\dots,T_d u).
\]
 Since
\[
\partial_t u = f + \sum^d_{p=1} A_pT_pu + A_0u,
\]
for any $j \in \{1,\dots,d\}$, applying the operator $T_j$ on the last
 equation we obtain
\[
\partial_t(T_ju) = \sum^d_{p=1} A_pT_p(T_ju)+ \sum^d_{p=1}
 (T_jA_p)T_pu + T_j(A_0u) + T_jf .
 \]
 We can then write
\[
\partial_tU = \sum^d_{p=1} B_p T_p U + B_0 U + F,
\]
 with
\[
F= (f, T_1 f,\dots,T_df),
\]
 and
\[
B_p = \left(\begin{array}{lllll}
 A_p &0 & \cdot &\cdot &0\\
 0 &A_p &0 &\cdot &\cdot\\
\cdot &0 &\cdot &\cdot &\cdot\\
\cdot &\cdot &\cdot &\cdot &0\\
 0 &\cdot &\cdot &0 &A_p
 \end{array}\right), \qquad p = 1,\dots,d,
 \]
 and the coef\/f\/icients of $B_0$ can be calculated from  the
 coef\/f\/icients of $A_p$ and from $T_jA_p$ with $(p = 0,\dots,d)$ and
 $(j=1,\dots,d)$.
 Using the induction hypothesis we then deduce  the result, and the
 proof of Lemma~\ref{lemma3.1} is f\/inished.
 \end{proof}

  {\bf B. Estimate about the approximated solution.}
  We notice that the necessary hypothesis to
 the proof of the inequalities of  Lemma~\ref{lemma3.1} require exactly one
 more derivative than the regularity which appears in the statement
 of the theorem that we have to prove. We then have  to regularize the system \eqref{S} by adapting the Friedrichs
 method. More precisely we consider the system 
\begin{gather}
 \left\{\begin{array}{l}
\displaystyle \partial_t
  u_n - \sum^d_{p=1} J_n(A_pT_p(J_nu_n)) - J_n(A_0 J_n
 u_n)  = J_nf,\\
  u{_n}|_{{t=0}}  = J_nu_0,
 \end{array}\right.\label{S_n}
 \end{gather}
 with $J_n$ is the cut of\/f operator  given by
\begin{gather*}J_n w = (J_n w_1,\dots,J_n w_m) \qquad  {\rm{and}} \qquad J_n w_j := {\cal F}^{-1}_D(1_{B(0,n)}(\xi)
  {\cal F}_D(w_j)), \qquad j=1,\dots,m.
  \end{gather*}
  Now we state the following proposition (see \cite[p.~389]{Chaz}) which we need in the sequel of this subsection.

  \begin{proposition}\label{proposition3.1}
  Let $E$ be a  Banach space, $I$ an open interval  of $\R$, $f \in C(I,E)$, $u_0 \in E$  and~$M$ be a continuous map from $I$ into $\mathcal{L}(E)$, the set of linear continuous applications from $E$ into itself. There exists a unique solution $u \in C^{1}(I,E)$  satisfying
  \begin{gather*}
   \left\{\begin{array}{l}
 \displaystyle  \frac{du}{dt}  = M(t)u + f,\vspace{1mm}\\
  u|_{{t=0}}  = u_0.
 \end{array}\right.
 \end{gather*}
  \end{proposition}

By taking $E = [L^2_k(
 \R^d)]^{m}$, and using the continuity of the operators $T_p J_n$  on $[L^2_k(
 \R^d)]^{m}$, we reduce the system \eqref{S_n} to  an evolution
  equation
\begin{gather*} \left\{\begin{array}{l}
 \displaystyle  \frac{du_n}{dt}  = M_n(t)u_n + J_n f,\vspace{1mm}\\
  {u_n}|_{{t=0}}  = J_nu_0
 \end{array}\right.
\end{gather*}
 on $[L^2_k(
 \R^d)]^{m}$, where
 \[
 t \mapsto M_n(t) = \displaystyle\sum^d_{p=1} J_n A_p(t,\cdot)T_pJ_n + J_n A_0(t,\cdot) J_n,
 \]
 is a continuous application from $I$ into   $\mathcal{L}([L^2_k(
 \R^d)]^{m})$. Then from Proposition~\ref{proposition3.1} there
 exists  a~unique function $u_n$ continuous on $I$ with values
 in $[L^2_k(
 \R^d)]^{m}$. Moreover, as the matrices $A_p$ are
 $C^{\infty}$ functions of $t$,
  $J_n f \in  [C(I, L^2_k(
 \R^d))]^m$ and $u_n$ satisfy
 \[
  \frac{du_n}{dt}  = M_n(t)u_n + J_n f.
\]
 Then $ \frac{du_n}{dt} \in [C(I, L^2_k(
 \R^d))]^m$ which implies that $u_n \in [C^1(I, L^2_k(
 \R^d))]^m$.
 Moreover, as $J^2_n =
 J_n$, it is obvious that $J_n u_n$ is also
 a solution of \eqref{S_n}. We apply Proposition~\ref{proposition3.1} we deduce that $J_nu_n = u_n$. The function $u_n$ is then belongs to $[C^1(I, H^s_k(
 \R^d))]^m$ for any integer $s$ and so  \eqref{S_n} can be written as
\begin{gather*}
\left\{ \begin{array}{l}
\displaystyle \partial_t u_n - \sum^d_{p=1} J_n(A_pT_pu_n) - J_n(A_0u_n)  =
 J_n f,\\
 u{_n}|_{{t=0}}  = J_n u_0.
 \end{array}\right.
 \end{gather*}
Now, let us estimate the evolution of $\|u_n(t)\|_{s,k}$.

 \begin{lemma} \label{lemma3.3} For any positive integer $s$, there exists a
 positive constant $\lambda_s$ such that for any integer $n$ and
 any $t$ in the interval $I$, we have
\[
\|u_n(t)\|_{s,k} \leq e^{\lambda_st}\|J_nu(0)\|_{s,k} +
 \int^t_0 e^{\lambda_s(t-t')}\|J_nf(t')\|_{s,k}dt' .
 \]
 \end{lemma}

 \begin{proof} The proof uses  the same ideas as in Lemma~\ref{lemma3.1}.
 \end{proof}

 {\bf C. Construction of solution.}
 This step consists on the proof of the following existence and  uniqueness  result:

\begin{proposition}\label{proposition3.2}
  For $s \geq 0$, we consider the
 symmetric system \eqref{S}
 with $f$ in $[C(I,\! H_k^{s+3}\!( \R^d))]^m\!\!$ and $v$  in $[H_k^{s+3}(
 \R^d)]^m$. There exists a unique solution $u$ belonging to the space
 $[C^1(I, H^s_k( \R^d))]^m\,\cap\, [C(I, H_k^{s+1}( \R^d))]^m$
 and satisfying the energy estimate
 \begin{gather} \|u(t) \|_{\sigma,k}
 \leq e^{\lambda_st}\|v\|_{\sigma,k} + \int^t_0 e^{\lambda_s(t-\tau)}
 \|f(\tau)\|_{\sigma,k}d\tau, \qquad {\rm{for \  all}}\ \ \sigma \leq s + 3 \ \  {\rm{and}}\ \ t \in I.\label{41}\end{gather}
 \end{proposition}

\begin{proof} Us consider the sequence $(u_n)_n$ def\/ined by the Friedrichs method and let us prove that this sequence is a Cauchy one in $[L^\infty(I, H^{s+1}_k( \R^d))]^m$. We put $V_{n,p} = u_{n+p} - u_n$,
 we have
\begin{gather*}
\left\{\begin{array}{l}
 \displaystyle \partial_t V_{n,p} - \sum^d_{j=1} J_{n+p} (A_jT_j V_{n,p}) -
 J_{n+p}(A_0 V_{n,p})  =  f_{n,p},\\
 {V_{n,p}}|_{t=0}  = (J_{n+p} - J_n)v
 \end{array}\right.
 \end{gather*}
 with
\[
f_{n,p} = - \displaystyle\sum^d_{j=1} (J_{n+p} - J_n)(A_jT_jV_{n,p}) - (J_{n+p} - J_n)(A_0 V_{n,p}
 )+(J_{n+p} - J_n)f.
 \]
 From Lemma~\ref{lemma3.3}, the sequence $(u_n)_{n \in \N}$ is
 bounded in $[L^\infty(I, H^{s+3}_k( \R^d))]^m$. Moreover, by
 a~simple calculation we f\/ind
\[
\|(J_{n+p} - J_n)(A_jT_jV_{n,p})\|_{s+1,k} \leq \frac{C}{n}
 \|A_jT_jV_{n,p}\|_{s+2,k} \leq \frac{C}{n}
  \|u_n(t)\|_{s+3,k}.
  \]
  Similarly, we have
\[
\|(J_{n+p} - J_n)(A_0V_{n,p}) + (J_{n+p} - J_n)f\|_{s+1,k} \leq \frac{C}{n}
  \big(\|u_n(t)\|_{s+3,k} + \|f(t)\|_{s+3,k}\big).
  \]
  By Lemma~\ref{lemma3.3} we deduce that
\[
\|V_{n,p} (t)\|_{s+1,k} \leq \frac{C}{n}  e^{\lambda_st} .
\]
  Then $(u_n)_n$ is  a Cauchy sequence in $[L^\infty(I, H_k^{s+1}
  ( \R^d))]^m$. We then have  the existence of a~solution $u$
  of \eqref{S}  in $[C(I, H_k^{s+1}( \R^d))]^m$. Moreover by the equation stated in~\eqref{S} we deduce that~$\partial_t u$ is in $[C(I, H_k^{s}( \R^d))]^m$, and so
  $u$ is in $[C^1(I,H_k^s(
   \R^d))]^m$. The uniqueness follows  immediately from Lemma~\ref{lemma3.3}.

   Finally we will  prove the inequality (\ref{41}).
   From Lemma~\ref{lemma3.3} we have
\[
\|u_n(t)\|_{s+3,k} \leq e^{\lambda_st}\|J_nu(0)\|_{s+3,k} +
   \int^t_0 e^{\lambda_s(t-\tau)}\|J_nf(\tau)\|_{s+3,k}d\tau.
   \]
  Thus
 \[
 \limsup_{n\rightarrow \infty} \|u_n(t)\|_{s+3,k} \leq e^{\lambda_st}
   \|v\|_{s+3,k} + \int^t_0 e^{\lambda_s(t-\tau)}\|f(\tau)\|_{s+3,k} d\tau .
 \]
   Since for any $t  \in I$, the sequence $(u_n(t))_{n \in \N}$ tends to $u(t)$  in $[H_k^{s+1}( \R^d)]^m$,  $(u_n(t))_{n \in
   \N}$ converge weakly to $u(t)$ in $[H_k^{s+3}(
   \R^d)]^m$, and then
\[
 u(t) \in [H_k^{s+3}( \R^d)]^m \mbox{ and }
 \|u(t)\|_{s+3,k} \leq \lim_{n\rightarrow \infty} \sup\|u_n(t)\|_{s+3,k} .
 \]
 The Proposition~\ref{proposition3.2} is thus proved.
 \end{proof}

 Now we will  prove the existence part of  Theorem~\ref{theorem3.1}.

 \begin{proposition}\label{proposition3.3}
 Let $s$ be an integer. If $v$ is in $[H^s_k( \R^d)]^m$
 and $f$ is in $[C(I,H^s_k( \R^d))]^m$, then there exists a
 solution of a symmetric system \eqref{S} in the space 
 \[
 [C(I,  H^s_k ( \R^d))]^m   \cap [C^1 (I,\!H^{s-1}_k ( \R^d))]^m.
 \]
 \end{proposition}

 \begin{proof} We consider the sequence $(\tilde{u}_n)_{n \in \N}$ of
 solutions of
\begin{gather*}
\left\{ \begin{array}{l}
\displaystyle \partial_t \tilde{u}_n - \sum^d_{j=1} (A_j T_j \tilde{u}_n) -
 (A_0 \tilde{u}_n)  = J_n f,\\
\tilde{u}_{n}|_{t=0}  = J_n v.
 \end{array}\right.
\end{gather*}
 From Proposition~\ref{proposition3.2} $(\tilde{u}_n)_n$ is in
  $[C^1(I, H_k^s(
 \R^d))]^m$.
 We will  prove that $(\tilde{u}_n)_n$ is a Cauchy sequence in $[L^\infty(I, H_k^s(
 \R^d))]^m$. We put $\tilde{V}_{n,p} = \tilde{u}_{n+p} -
 \tilde{u}_n$. By dif\/ference, we f\/ind
 \begin{gather*}
 \left\{ \begin{array}{l}
\displaystyle \partial_t \tilde{V}_{n,p} - \sum^d_{j=1} A_j T_j
 \tilde{V}_{n,p} - A_0\tilde{V}_{n,p}  = (J_{n+p} - J_n )f,\\
 {\widetilde{V_{n,p}}}|_{t=0} = (J_{n+p} - J_n)v .
 \end{array}\right.
 \end{gather*}
 By Lemma~\ref{lemma3.3} we deduce that
 \[
 \|\tilde{V}_{n,p}\|_{s,k} \leq e^{\lambda_st}
 \|(J_{n+p} - J_n)v\|_{s,k} + \int^t_0 e^{\lambda_s(t-\tau)}
 \|(J_{n+p} - J_n)f(\tau)\|_{s,k}d\tau .
 \]
Since $f$ is in $[C(I, H^s_k( \R^d))]^m$, the sequence
 $(J_nf)_n$ converges to $f$ in $[L^\infty([0,T], H^s_k(
 \R^d))]^m$, and since $v$ is in $[H^s_k(
 \R^d)]^m$, the sequence $(J_nv)_n$ converge to $v$ in
 $[H^s_k(\R^d)]^m$ and so
  $(\tilde{u}_n)_n$ is a Cauchy sequence in $[L^\infty(I, H^s_k(
 \R^d))]^m$. Hence it converges to  a function $u$ of $[C(I, H^s_k(
 \R^d))]^m$, solution of the system~\eqref{S}. Thus $\partial_tu$
 is in $[C(I, H_k^{s-1}( \R^d))]^m$ and the proposition is
proved.
\end{proof}

 The existence in  Theorem~\ref{theorem3.1} is then
proved as well as the uniqueness, when $s \geq 1$.

 {\bf D. Uniqueness of solutions.}
 In the following we give the result of uniqueness for $s = 0$ and
 hence  Theorem~\ref{theorem3.1} is proved.

 \begin{proposition} \label{proposition3.4} Let $u$ be a solution in $[C(I, L^2_k(
 \R^d))]^m$   of the symmetric system
\begin{gather*}
\left\{ \begin{array}{l}
\displaystyle \partial_t u - \sum^d_{j=1} A_jT_ju - A_0u  = 0,\\
 u|_{ t=0} = 0.
 \end{array}\right.
 \end{gather*}
 Then $u \equiv 0$.\end{proposition}

 \begin{proof} Let $\psi$ be a function in $[D(]0,T[ \times \R^d)]^m$; we
 consider the following system
\begin{gather}
\left\{ \begin{array}{ll}
\displaystyle -\partial_t\varphi + \sum^d_{j=1} T_j(A_j\varphi) - {}^{t}A_{0}\varphi  = \psi,\\
 \varphi|_{t=T}  = 0.
 \end{array}\right.\label{{}^tS}
 \end{gather}
 Since
\[
T_j(A_j\varphi) = A_jT_j\varphi + (T_jA_j) \varphi,
\]
  the system \eqref{{}^tS}  can be written
\begin{gather}
\left\{ \begin{array}{ll}
\displaystyle - \partial_t\varphi + \sum^d_{j=1} A_jT_j\varphi -
\tilde{A}_0\varphi  =  \psi,\\
\varphi|_{t=T}  = 0
 \end{array}\right.\label{{}^tS'}
 \end{gather}
 with
\[
\tilde{A}_0 =  {}^tA_0 - \sum^d_{j=1} T_jA_j .
\]
 Due to  Proposition \ref{proposition3.2}, for any integer $s$ there exists a solution $\varphi$ of
\eqref{{}^tS'}    in
  $[C^1([0,T],\! H^s_k(
 \R))]^m$.
We  then have
 \begin{gather*}
 \langle u,\psi \rangle_k  =  \langle u,
  -\partial_t\varphi +
 \sum^d_{j=1} A_jT_j\varphi - \tilde{A}_0\varphi\rangle_k\\
\phantom{\langle u,\psi \rangle_k }{} = - \int_I \langle u(t,\cdot), \partial_t\varphi(t,\cdot)\rangle_k dt
 +
 \sum^d_{j=1} \int_{I \times \R^d}
 u(t,x)T_j(A_j\varphi)(t,x) \omega_k(x)dtdx\\
\phantom{\langle u,\psi \rangle_k  =}{} -\int_{I \times \R^d} u(t,x)\,{}^{t} A_0\varphi(t,x)
 \omega_k(x)dt dx
 \end{gather*}
  with $\langle\cdot,\cdot\rangle_k$ def\/ined by
\[
\langle u,\chi\rangle_k = \dint_{I}\langle u(t,\cdot),\chi(t,\cdot)\rangle_k dt =  \dint_{I\times\R^d}u({t},x)\chi(t,x)\omega_k(x)dx dt,
\qquad \chi \in {\cal S}(\R^{d+1}).
\]
 By using that $u(t,\cdot)$ is in $[L^2_k( \R^d)]^m$ for any
 $t$ in $I$ and the fact that $A_j$ is symmetric we obtain
\[
\int_{I \times \R^d} u(t,x)T_j(A_j\varphi)(t,x)
 \omega_k(x)dt dx =
  - \int_I \langle A_jT_j u(t,\cdot), \varphi(t,\cdot)
  \rangle_k dt.
  \]
  So
\[
\langle u, \psi\rangle _k = - \int_I \langle u(t,\cdot ),\partial_t\varphi(t,\cdot)
  \rangle_k dt - \sum^d_{j=1} \langle A_jT_j u + A_0 u,\varphi\rangle_k.
\]
  As $u$ is not very regular, we have to justify
   the integration by parts
  in time on the quantity
  $\dint_I \langle u(t,\cdot), \partial_t \varphi(t,\cdot) \rangle_k dt.$
  Since  $J_n u(\cdot,x)$  are
   $C^1$ functions on
  $I$, then by integration by parts, we obtain, for any $x \in
  \R^d$,
\begin{gather*}
- \int_I J_nu(t,x)\partial_t\varphi(t,x)dt = - J_nu(T,x)\varphi(T,x) + J_nu(0,x) \varphi(0,x)
  + \int_I \partial_t J_nu(t,x) \varphi(t,x)dt .
  \end{gather*}
  Since $u(0,\cdot) = \varphi(T,\cdot) = 0$, we have
\[
- \int_I J_nu(t,x)\partial_t\varphi(t,x) dt
   = \int_I\partial_t(J_nu)(t,x)\varphi(t,x)dt .
   \]
   Integrating with respect to
$\displaystyle\omega_{k}(x)dx$ we obtain
   \begin{equation}
   - \int_{I \times \R^d} J_n u(t,x) \partial_t \varphi(t,x) \omega_k(x)dtdx = \int_I \langle \partial_t(J_nu)(t,\cdot)
   ,\varphi(t,\cdot)\rangle_k dt.\label{42}
   \end{equation}
   Since $u$ is in $[C(I, L^2_k( \R^d))]^m \cap [C^1(I, H_k^{-1}(
   \R^d))]^m$, we have
\[
\lim_{n\rightarrow \infty} J_n u = u \mbox{ in }
   [L^\infty(I, L^2_k( \R^d))]^m \qquad \mbox{and}\qquad
   \lim_{n\rightarrow \infty} J_n\partial_tu = \partial_t u \quad \mbox{in} \ \
   [L^\infty(I, H_k^{-1}( \R^d))]^m.
   \]
   By passing  to the limit in (\ref{42}) we obtain
\[
- \int_I \langle u(t,\cdot), \partial_t \varphi(t,\cdot)\rangle_k dt = \int_I
   \langle \partial_t u(t,\cdot),\varphi(t,\cdot)\rangle_kdt .
   \]
   Hence
\[
\langle u,\psi\rangle_k = \int_I \langle \partial_t u(t,\cdot) -
   \sum^d_{j=1} \langle A_j T_j u(t,\cdot) - A_0 u(t,\cdot), \varphi(t,\cdot)\rangle_k dt.
\]
   However since $u$ is a solution of \eqref{S} with $f \equiv 0$, then $u \equiv
   0$. This ends  the proof.\end{proof}

   \subsection[The Dunkl-wave equations with variable coefficients]{The Dunkl-wave equations with variable coef\/f\/icients}\label{sec3.2}

   For $t \in \R$ and $x \in \R^d$, let $P(t, x, \partial_t,
   T_x)$ be a dif\/ferential-dif\/ference operator of degree
    2 def\/ined by
   \begin{equation*}P(u) = \partial^2_t u - {\rm{div}}_{k}[A\cdot \nabla_{k,x}u] + Q(t,x,\partial_tu, T_x u),
   \end{equation*}
   where
\[
\nabla_{k,x}u := \left(T_1  u,\dots,
T_d\, u\right), \qquad {\rm{div}}_{k}  \left(  v_1,\dots,v_d  \right) :=
 \sum_{i=1}^{d}T_i v_i,
 \]
    $A$ is a real symmetric matrix such that there exists $m > 0$ satisfying
   \begin{equation*}
   \langle A(t,x)\xi, \xi \rangle \geq m\|\xi\|^2,\qquad {\rm{for \  all}}\ \ (t,x)\in \R\times \R^{d}  \ \ {\rm{and}} \ \ \xi \in \R^d 
   \end{equation*}
   and $Q(t,x,\partial_tu, T_x u)$ is dif\/ferential-dif\/ference operator of degree~1, and the
   matrix $A$ is  $W$-invariant with respect to~$x$;  the coef\/f\/icients of
   $A$ and $Q$ are  $C^{\infty}$ and all derivatives are bounded. If we
   put $B = \sqrt{A}$ it is easy to see that
  the
   coef\/f\/icients of $B$ are $C^\infty$ and all derivatives are  bounded.

   We introduce the vector $U$ with $d+2$ components
   \begin{equation*} U = \left(
   u,
   \partial_t u,
   B\nabla_{k,x}u
   \right).
   \end{equation*}
   Then, the equation $P(u) = f$ can  be written as
   \begin{equation}\partial_tU = \left(\sum^d_{p=1} A_p T_p\right) U + A_0 U
    + \left(
   0,f,0
   \right)\label{46}\end{equation}
   with
\[
A_p = \left(\begin{array}{lllllll}
   0 &\cdot  &\cdot &\cdot &\cdot &\cdot &0\\
   \cdot&0 &b_{p\, 1} &\cdot &\cdot &\cdot &b_{p\, d}\\
   \cdot &b_{1 \,p} &0 &\cdot &\cdot &\cdot &0\\
   \cdot &\cdot &\cdot &\cdot &\cdot &\cdot &\cdot\\
   \cdot &\cdot &\cdot &\cdot &\cdot &\cdot &\cdot\\
   0 &b_{d\, p} &0 &\cdot &\cdot &\cdot &0
   \end{array}\right)
   \]
   and $B = (b_{ij})$. Thus the system (\ref{46}) is symmetric and  from Theorem~\ref{theorem3.1} we deduce the following.

  \begin{theorem}\label{theorem3.2} For all $s \in \N$ and $u_0 \in H_k^{s+1} ( \R^d),\;
   u_1 \in H^s_k( \R^d)$ and $f \in C( \R, H^s_k(
   \R^d))$, there exists a unique $u \in C^1( \R, H^s_k( \R^d)) \cap
   \,C( \R, H_k^{s+1} ( \R^d))$ such that
\begin{gather*}
\left\{ \begin{array}{l}
   \partial^2_t u - {\rm{div}}_{k}[A\cdot \nabla_{k,x} u] + Q(t,x,\partial_t u,
   T_xu)  = f,\\
   u|_{t=0}  =u_0,\\
   \partial_t u|_{t = 0} = u_1.
   \end{array}\right.
\end{gather*}
   \end{theorem}

   \subsection{Finite speed of propagation}\label{sec3.3}

   \begin{theorem}\label{theorem3.3} Let \eqref{S} be a symmetric system. There exists a
   positive constant $C_0$ such that, for any positive real $R$,
    any function $f \in [C(I, L^2_k( \R^d))]^m$ and any $v \in [L^2_k(
   \R^d)]^m$ satisfying
    \begin{gather}
   f(t,x) \equiv  0 \qquad {\rm{for}}\quad \|x\| < R - C_0t, \label{hy1} \\ \label{hy2}
   v(x) \equiv  0 \qquad{\rm{for}}\quad \|x\| < R,
   \end{gather}
   the unique solution $u$ of system \eqref{S} belongs to $[C(I, L^2_k(
   \R^d))]^m$ with
\[
u(t,x) \equiv 0\qquad {\rm{for}}\quad \|x\| < R - C_0 t.
\]
   \end{theorem}

  \begin{proof} If $f_\varepsilon \in [C(I, H^1_k( \R^d))]^m$, $v_\varepsilon \in [H^1_k( \R^d))]^m$ are given such that $f_\varepsilon \to f$ in $[C(I, L^2_k( \R^d))]^m$ and $v_\varepsilon \!\to\! v$ in $[L^2_k(
   \R^d\!)]^m$, we know by Subsection~\ref{sec3.1} that
    the solution $u_\varepsilon$ belongs to  $[C(I,\! H^1_k\!(
   \R^d\!))]^m\!$ and satisf\/ies  $u_\varepsilon \to u$ in $[C(I, L^2_k( \R^d))]^m$. Therefore, if we construct $f_\varepsilon$ and $v_\varepsilon$ satisfying
   (\ref{hy1}) and (\ref{hy2}) with $R$ replaced by $R-\varepsilon$, it suf\/f\/ices to prove Theorem~\ref{theorem3.3} for $f \in [C(I, H^1_k( \R^d))]^m$ and $v \in [H^1_k( \R^d))]^m$. We have then $u \in [C^{1}(I, L^2_k( \R^d))]^m \bigcap [C(I, H^1_k( \R^d))]^m$. To this end let us consider
    $\chi \in D(\R^{d})$ radial with ${\rm supp}\,\chi \subset  B(0,1)$ and
\[
 \int_{ \R^d} \chi(x) \omega_k(x)dx = 1.
 \]
   For $\varepsilon > 0$, we put
\begin{gather*}
u_{0,\varepsilon} = \chi_\varepsilon\ast_D v := (\chi_\varepsilon\ast_D v_1,\dots,\chi_\varepsilon\ast_D v_d),\\
 f_\varepsilon(t,\cdot) = \chi_\varepsilon \ast_Df(t,\cdot) := (\chi_\varepsilon \ast_D f_1(t,\cdot),\dots,\chi_\varepsilon \ast_D f_d(t,\cdot)),
 \end{gather*}
with
\[
\chi_\varepsilon(x) = \frac{1}{\varepsilon^{d+2\gamma}}
   \chi\left(\frac{x}{\varepsilon}\right).
\]
   The hypothesis (\ref{hy1}) and (\ref{hy2}) are then satisf\/ied by
    $f_\varepsilon$ and $u_{0,\varepsilon}$ if we replace $R$ by $R-\varepsilon$.
   On the other hand the solution $u_\varepsilon$ associated with
   $f_\varepsilon$ and $u_{0,\varepsilon}$   is $[C^1(I, H^s_k(
   \R^d))]^m$ for any integer~$s$. For $\tau \geq 1$, we put
\[
u_\tau(t,x) = \exp\big(\tau(-t+\psi(x))\big)u(t,x),
\]
   where the function $\psi \in C^{\infty}( \R^d)$
    will be chosen later.

   By a simple calculation we see that
\[
\partial_tu_\tau - \sum^d_{j=1} A_jT_ju_\tau - B_\tau u_\tau = f_\tau
\]
   with
\[
f_\tau(t,x) = \exp(\tau(-t+\psi(x))) f(t,x),\qquad  B_\tau = A_0 + \tau\left(-{\rm Id} - \sum^d_{j=1} (T_j \psi)A_j\right).
\]
   There exists a positive constant $K$ such  that if $\|T_j\psi\|_{L^\infty_k
   ( \R^d)} \leq K$ for any $j=1,\dots,d$, we have for any $(t,x)$
\[
\langle {\rm{Re}}(B_\tau y), \bar{y}\rangle \leq \langle {\rm{Re}}(A_0 y), \bar{y}\rangle\qquad
{\rm{for \  all}} \ \ \tau \geq 1 \ \  {\rm{and}} \ \ y \in \C^m .
\] We proceed as
in the proof of energy estimate (\ref{1+}), we obtain the existence of
positive constant~$\delta_0$, independent of $\tau$, such that
for any $t$ in $I$, we have
\begin{equation}\|u_\tau(t)\|_{0,k} \leq e^{\delta_0 t} \|u_\tau(0)\|_{0,k}
+ \int^t_0 e^{\delta_0(t-t')}\|f_\tau(t')\|_{0,k}dt' .\label{47}\end{equation}
We put $C_0 = \frac{1}{K}$ and choose $\psi = \psi(\|x\|)$ such
that $\psi$ is  $C^\infty$ and such that
\[
-2\varepsilon + K(R - \|x\|) \leq \psi(x)
\leq - \varepsilon
 + K(R - \|x\|).
 \]
 There exists $\varepsilon > 0$ such that $\psi(x) \leq - \varepsilon + K(R -
 \|x\|)$. Hence
\[
 \|x\| \geq R - C_0 t, \qquad {\rm{for \  all}} \ \  (t,x) \quad \Longrightarrow \quad  -t +
 \psi(x) \leq - \varepsilon .
 \]
 Let   $\tau$ tend to $+ \infty$ in (\ref{47}), we deduce that
\[
 \lim_{\tau + \infty} \int_{\R^d}
 \exp(2\tau(-t + \psi(x)))||u(t,x)||^2 \omega_k(x)dx = 0, \qquad {\rm{for \  all}} \ \  t \in I.
\]
 Then
\[
u(t,x) = 0\qquad \mbox{on} \ \ \big\{(t,x) \in I \times \R^d ;\; t < \psi(x)\big\}.
\]
 However  if $(t_0,x_0)$ satisf\/ies  $\|x_0\| < R - C_0t_0$, we can f\/ind a
 function $\psi$ of previous type such that $t_0 < \psi(x_0)$.
 Thus the theorem is proved.
 \end{proof}

 \begin{theorem}\label{theorem3.4} Let \eqref{S} be a symmetric system. We assume
 that the functions $f {\in} [C(I,\! L^2_k\!( \R^d\!))]^m\!\!$ and $v \in [L^2_k(
 \R^d)]^m$ satisfy
\begin{gather*}
 f(t,x) \equiv  0 \qquad \mbox{for} \ \  \|x\| > R + C_0t,\\
 v(x)  \equiv  0 \qquad \mbox{for} \ \ \|x\| > R.
 \end{gather*}
 Then the unique solution $u$ of system \eqref{S} belongs to $[C(I, L^2_k(
 \R^d))]^m$ with
\[
u(t,x) \equiv 0\qquad \mbox{for}\ \ \|x\| > R + C_0 t.
\]
 \end{theorem}

 \begin{proof} The proof uses  the same ideas as in Theorem~\ref{theorem3.3}.\end{proof}

 \section{Semi-linear Dunkl-wave equations}\label{sec4}

 We consider the problem \eqref{P}. 
We denote by
 $ \|\Lambda_{k}u(t,\cdot)\|_{L_k^\infty}$ the norm def\/ined by
\[
  \|\Lambda_{k}u(t,\cdot)\|_{L_k^\infty} =  \|\partial_{t}u(t,\cdot)\|_{L_k^\infty(\R^{d})} + \sum_{j=1}^{d}\|T_{j}u(t,\cdot)\|_{L_k^\infty(\R^{d})}.
\]
 $ \|\Lambda_{k}u(t,\cdot)\|_{s,k}$ the norm def\/ined by
\[
  \|\Lambda_{k}u(t,\cdot)\|_{s,k}^{2} =  \|\partial_{t}u(t,\cdot)\|_{H_k^s(\R^{d})}^{2} + \sum_{j=1}^{d}\|T_{j}u(t,\cdot)\|_{H_k^s(\R^{d})}^{2}.
\]
 The main  result of this section is the following:
 \begin{theorem}\label{theorem4.1} Let $(u_0, u_1)$ be in $H^s_k( \R^d)
 \times H^{s-1}_{k}( \R^d)$ for $s > \gamma + \frac{d}{2}
 +1$. Then there exists a~positive time $T$  such that the problem \eqref{P} has a unique solution $u$ belonging to
\[
C([0,T] , H^s_k( \R^d)) \cap C^1([0,T] , H^{s-1}_k( \R^d))
\]
  and satisfying  the following blow up criteria: if $T^{*}$ denotes the maximal time of existence of such a solution, we have:

   -- The existence of  constant $C$ depending  only
 on $\gamma$, $d$ and  quadratic form $Q$ such that
 \begin{equation}
 T^{*} \geq \frac{C(\gamma, d,Q)}{\|\Lambda_k u(0,\cdot)\|_{{s-1},k}}.
 \label{48}\end{equation}

-- If $T^{*} < \infty$, then
 \begin{equation} \dint^{T^\ast}_0 \|\Lambda_{k}u(t,\cdot)\|_{L_k^\infty}
 dt = + \infty .\label{49}\end{equation}
 \end{theorem}

To prove  Theorem~\ref{theorem4.1} we need the following lemmas.

 \begin{lemma}\label{lemma4.1}
  {\rm (Energy Estimate in $H^k_s( \R^d)$.)}
If $u$ belongs to $C^1( I, H^s_k( \R^d)) \cap C(I, H_k^{s+1}(
 \R^{d}))$ for an integer $s$ and with $f$ defined by
\[
f = \partial^2_t u - \Delta_k u
\]
 then we have
 \begin{equation}\|\Lambda_{k}u(t,\cdot)\|_{s-1,k} \leq \|\Lambda_{k}u(0,\cdot)\|_{s-1,k} + \int^t_0 \|f(t',\cdot)\|_{H_k^{s-1}( \R^d)} dt', \qquad \mbox{for} \ \  t \in I.
 \label{50}\end{equation}
 \end{lemma}

 \begin{proof} We multiply the equation by $\partial_t u$ and we obtain
\[
 (\partial^2_t u, \partial_t u)_{H_k^{s-1}( \R^d)}
 - (\Delta_k u, \partial_t u)_{H_k^{s-1}( \R^d)}
 = \langle f, \partial_t u\rangle_{H_k^{s-1}( \R^d)}.
 \]
 A simple calculation yields that
\[
- \langle \Delta_k u, \partial_t u\rangle_{H_k^{s-1}( \R^d)}
 = \langle \nabla_k u, \nabla_k \partial_t u\rangle_{H_k^{s-1}
 ( \R^d)} .
 \]
 Thus
\[
\frac{1}{2} \frac{d}{dt} \|\Lambda_{k}u\|^2_{{s-1},k}
 = \langle f, \partial_t u\rangle_{H_k^{s-1}(\R^d)}.
 \]
 If $f = 0$, we deduce the conservation of energy
\[
\|\Lambda_{k}u(t,\cdot)\|^2_{{s-1},k}
  = \|\Lambda_{k}u(0,\cdot)\|^2_{{s-1},k}.
\]
  Otherwise,  Lemma~\ref{lemma3.2} gives
  \begin{equation*}\|\Lambda_{k}u(t,\cdot)\|_{{s-1},k}
  \leq \|\Lambda_{k}u(0,\cdot)\|_{{s-1},k} + \dint^t_0 \|f(t',\cdot)\|_{H_k^{s-1}
  ( \R^d)}dt'.\tag*{\qed} 
  \end{equation*}\renewcommand{\qed}{}
  \end{proof}

\begin{lemma}\label{lemma4.2} Let $(u_n)_{n\in \N}$ be the sequence
 defined by
\begin{gather*}
 \partial^2_t u_{n+1} - \Delta_k u_{n+1} = Q(\Lambda_{k}u_n, \Lambda_{k}u_n),\\
 (u_{n+1}, \partial_t u_{n+1})|_{t=0} = (S_{n+1} u_0, S_{n+1}u_1),
 \end{gather*}
 where $u_0 = 0$ and $S_{n+1}u_j$ defined by
\[
{\cal F}_D(S_{n+1}u_j)(\xi) = \psi(2^{-(n+1)}\xi)
 {\cal F}_D(u_j)(\xi) ,
\]
 with $\psi$ a function of $D( \R^d)$ such that $0 \leq \psi \leq
 1$ and ${\rm{supp}}\,\psi \subset B(0,1)$.

   Then there exists a positive time $T$  such that, the sequence $(\Lambda_{k}u_n)_{n \in
  \N}$ is bounded in the space $[L^\infty([0,T], H_k^{s-1}(
  \R^d))]^{d+1}$.
  \end{lemma}

  \begin{proof} First, from  Theorem~\ref{theorem3.2} the sequence $(u_n)_n$ is well def\/ined.  Moreover,
 due to the energy estimate,
\[
\|\Lambda_{k}u_{n+1}(t,\cdot)\|_{{s-1},k}
  \leq \|S_{n+1} \theta\|_{{s-1},k} + \dint^t_0
  \|Q(\Lambda_{k}u_n(t',\cdot ),  \Lambda_{k}u_n(t',\cdot))\|_{H_k^{s-1}( \R^d)}dt',
\]
  where
\[
\theta = (u_1,\nabla_k u_0) = \Lambda_{k}u(0,\cdot).
\]
  Since $s-1 > \gamma + \frac{d}{2}$, then from Theorem~\ref{theorem2.2}\,ii, we
  have
  \begin{equation}\|\Lambda_{k}u_{n+1}(t,\cdot)\|_{{s-1},k} \leq
  \|\theta\|_{{s-1},k} + C \int^t_0 \|\Lambda_{k}u_n(\tau,\cdot)\|^2_{{s-1},k}
  d\tau .\label{54}
  \end{equation}
  Let $T$ be a positive real such that
 \begin{equation}
 4CT \|\theta\|_{{s-1},k} < 1.\label{55}
 \end{equation}
  We will  prove by induction that, for any integer $n$
   \begin{equation}
   \|\Lambda_{k}u_{n+1}(t,\cdot)\|_{{s-1},k}
  \leq 2\|\theta\|_{{s-1},k} .\label{56}\end{equation}
  This property  is true for $n = 0$. We assume that it is true for~$n$.
  With the inequalities~(\ref{54}) and~(\ref{55}), we deduce that, for all $t \leq
  T$, we have
\begin{gather*}
  \|\Lambda_{k}u_{n+1}(t,\cdot)\|_{{s-1},k}  \leq  \|\theta\|_{{s-1},k}+ 4CT\|\theta\|_{{s-1},k}^{2} \leq (1+4CT\|\theta\|_{
  {s-1},k})\|\theta\|_{{s-1},k},\\
  \|\Lambda_{k}u_{n+1}(t,\cdot)\|_{H_k^{s-1}( \R^d)}  \leq  2   \|\theta\|_{{s-1},k}.
  \end{gather*}
  This gives (\ref{56})
  and the proof of lemma is established.
\end{proof}

\begin{lemma} \label{lemma4.3} There exists a positive time $T$
 such that,
  $(\Lambda_{k}u_n)_{n \in
  \N}$ is a Cauchy sequence in the space $[L^\infty([0,T], H_k^{s-1}(
  \R^d))]^{d+1}$.
  \end{lemma}
  \begin{proof} We put
\[
V_{n,p} = u_{n+p} - u_n .
\]
  By dif\/ference, we see that
\begin{gather*}
\left\{\begin{array}{l}  \partial^2_t V_{n+1,p} - \Delta_k V_{n+1,p}  = Q(\Lambda_{k}V_{n,p},
  \Lambda_{k}u_{n+p}+ \Lambda_{k}u_n),\vspace{1mm}\\
  {\Lambda_{k}V_{n+1,p}}|_{t=0} = (S_{n+p+1} - S_{n+1})\theta.
  \end{array}\right.
  \end{gather*}
  By energy estimate, we establish from  (\ref{56}) that
   \begin{gather*}
   \|\Lambda_{k}V_{n+1,p}(t,\cdot)\|_{{s-1},k} \leq \|(S_{n+p+1} -
  S_{n+1})\theta\|_{{s-1},k} \nonumber\\
  \phantom{\|\Lambda_{k}V_{n+1,p}(t,\cdot)\|_{{s-1},k} \leq}{}
  + 4CT\|\theta\|_{{s-1},k}
  \|\Lambda_{k}V_{n,p}\|_{[L^\infty([0,T],H_k^{s-1}( \R^d))]^{d+1}}.
  \end{gather*}
  We put
\[
\rho_n = \sup_{p \in \N}\|\Lambda_{k}V_{n,p}\|_{[L^\infty([0,T],H_k^{s-1}( \R^d))]^{d+1}}
   \qquad {\rm{and}} \qquad \varepsilon_{n} = \sup_{p \in \N} \|(S_{n+p} - S_n)\theta\|_{{s-1},k}.
\]
  From this and the last inequality, we have
\[
\rho_{n+1} \leq \varepsilon_{n+1} + 4CT \rho_n\|\theta\|_{{s-1},k}.
\]
  The sequence $(S_n\theta)_{n \in \N}$ converges to
  $\theta$ in $[H_s^{s-1}( \R^d)]^{d+1}$.
  By passing to the superior
  limit, we obtain
\[
\limsup_{n\rightarrow \infty} (\rho_{n+1})
  \leq 0 + 4CT\|\theta\|_{{s-1},k}
  \limsup_{n\rightarrow \infty} (\rho_n).
\]
  However, since
\[
\limsup_{n\rightarrow \infty}(\rho_{n+1}) = \limsup_{n\rightarrow \infty} (\rho_n),
\]
  we deduce
\[
\limsup_{n\rightarrow \infty} (\rho_n) \leq 4CT\|\theta\|_{{s-1},k}
  \limsup_{n\rightarrow \infty} (\rho_n),
  \]
  and the result holds by $4CT\|\theta\|_{{s-1},k} <
  1$.

  Hence
\[
\limsup_{n\rightarrow \infty}(\rho_n) = 0.
\]
  Then  $(\Lambda_{k}u_n)$ is a Cauchy sequence in $[L^\infty([0,T], H_k^{s-1}
  ( \R^d))]^{d+1}$, and so  Lemma~\ref{lemma4.3} is proved.
  \end{proof}

\begin{proof}[Proof of Theorem~\ref{theorem4.1}]
  In the following we will  prove that the unique solution  of~\eqref{P} belongs to
   $C([0,T],H_k^{s-1}\!(\R^d))$. Indeed, Lemma~\ref{lemma4.3} implies the existence of $v$ in $[L^\infty([0,T], H_k^{s-1}(\R^d))]^{d+1}$ such that
\[
\Lambda_{k}u_n \to v \qquad \mbox{in}\ \ [L^\infty([0,T], H_k^{s-1}(\R^d))]^{d+1}.
\]
Moreover,  Lemma~\ref{lemma4.2} gives the existence of a positive time $T$ such that the sequence $(u_n)_n$ is bounded in
$L^\infty([0,T], H_k^{s}(\R^d))$. Thus there exists $u$ such that the sequence $(u_n)_n$ converges weakly  to $u$ in $L^\infty([0,T], H_k^{s}(\R^d))$.

 The uniqueness for the solution of \eqref{P} gives that $v = \Lambda_{k}u$ and that
 \[
 u_n \to u \qquad \mbox{in}\ \  L^\infty([0,T], H_k^{s}(\R^d)).
 \]
Finally it is easy to see that $u$ is the unique solution of~\eqref{P} which   belongs to $C([0,T], H_k^{s}(\R^d))$.

Now we are going to prove the inequalities (\ref{48}).
  We have proved that if
  $T < \frac{1}{4C\|\theta\|_{{s-1},k}}$, then $u \in C([0,T], H_k^{s}(\R^d))$. Hence, if $T^{*}$ denote the maximal time of existence of such a solution we have $T^{*} > T$ and this gives that $T^{*} \geq \frac{1}{4C\|\theta\|_{
  {s-1},k}}$
  and $u \in C([0,T^{*}[, H_k^{s}(\R^d))$.
   Finally we will  prove the condition (\ref{49}). We assume that $T^{*} < \infty$ and $\dint_{0}^{T^{*}}
   \|\Lambda_k u(t)\|_{
   {s-1},k}dt < \infty$.
   Indeed, it is easy to see that the
   maximal time of   solution of problem \eqref{P}  with initial
   given $u(t)$ is $T^\ast -t$.
   Thus, from the relation (\ref{48}) we deduce that
\[
T^\ast - t \geq \frac{C}{\|\Lambda_k u(t,\cdot )\|_{s-1,k}} .
\]
   This implies that
   \begin{equation} \|\Lambda_k u(t,\cdot)\|_{s-1,k} \geq \frac{C}{T^\ast -t},\qquad {\rm{for \  all}} \ \ t \in [0,T^{*}[.\label{59}\end{equation}
   Hence $\|\Lambda_k u(t,\cdot)\|_{s-1,k}$ is not bounded if $t$
   tends to $T^\ast$.

   On the other hand from (\ref{50}) and Theorem~\ref{theorem2.2}\,i there exists  a positive constant $C$
    such that
\begin{gather*} \|\Lambda_k u(t)\|_{{s-1},k}
    \leq \|\Lambda_k u(0)\|_{
    {s-1},k} \\
    \phantom{ \|\Lambda_k u(t)\|_{{s-1},k}\leq}{} + C\dint^t_0 \|\Lambda_k u(t',\cdot)\|_{L^\infty_k}\|\Lambda_k u(t',\cdot)\|_{{s-1},k} dt',
    \qquad {\rm{for \ all}} \ \ t \in [0,T^{*}[.
\end{gather*}
    Then from the usual Gronwall lemma we obtain
    \begin{gather}\|\Lambda_k u(t,\cdot)\|_{{s-1},k}
    \leq C\|\Lambda_k u(0,\cdot)\|_{{s-1},k}\nonumber\\
    \phantom{\|\Lambda_k u(t,\cdot)\|_{{s-1},k}\leq}{}\times{} \exp\left(\dint^t_0 \|\Lambda_k u(t',\cdot)\|_{L^\infty_k
    } dt'\right),\qquad {\rm{for \ all}} \ \ t \in [0,T^{*}[.\label{las}\end{gather}
  Finally if we tend  $t$ to $T^*$ in (\ref{las}) we obtain that $\|\Lambda_k u(t)\|_{{s-1},k}$ is  bounded which is not true from~(\ref{59}).
  Thus we have proved~(\ref{49}) and the proof of Theorem~\ref{theorem4.1} is f\/inished.
  \end{proof}

 \subsection*{Acknowledgements} We would like to thank Professor
 Khalifa Trim\`eche for his help and encouragement. Thanks are also due to the referees and editors for their suggestions and comments.

\pdfbookmark[1]{References}{ref}
\LastPageEnding
\end{document}